\magnification 1200
\def\R{{\rm I\kern-0.2em R\kern0.2em \kern-0.2em}}
\def\N{{\rm I\kern-0.2em N\kern0.2em \kern-0.2em}}
\def\P{{\rm I\kern-0.2em P\kern0.2em \kern-0.2em}}
\def\B{{\rm I\kern-0.2em B\kern0.2em \kern-0.2em}}
\def\Z{{\rm I\kern-0.2em Z\kern0.2em \kern-0.2em}}
\def\C{{\bf \rm C}\kern-.4em {\vrule height1.4ex width.08em depth-.04ex}\;}

\def\D{{\Delta}}

\def\z{{\zeta}}

\def\cW{{\cal W}}

\font\ninerm=cmr8
\ 
\vskip 22mm
\centerline {\bf ON MEROMORPHIC EXTENDIBILITY}
\vskip 8mm
\centerline{Josip Globevnik}
\vskip 8mm
{\noindent \ninerm ABSTRACT \ \ Let $D$ be a bounded domain in the complex plane 
whose boundary consists of finitely many 
 pairwise disjoint real analytic simple closed curves. Let $f$ be an integrable function on $bD$. In the 
 paper we show how to compute the candidates for poles of a meromorphic extension of $f$ through $D$ and 
 thus reduce the question of meromorphic extendibility to the question of holomorphic extendibility. 
Let 
$A(D)$ be the algebra of all continuous 
functions on $\overline D$ which are holomorphic on $D$. We prove that a 
continuous function $f$ on $bD$ extends meromorphically through 
$D $ if and only if there is an $N\in\N\cup \{0 \} $ such that 
the change of argument of $Pf+Q$ along $bD$ is bounded below by $-2\pi N$  
for all  $P, Q \in A(D) $ such that $Pf+Q\not= 0 $ on $bD $. If this 
is the case then the meromorphic extension of $f$ has at most $N$ poles 
in $D$, counting multiplicity. }
\vskip 6mm
\bf 1.\ Introduction \rm 
\vskip 2mm
Let $D\subset \C$ be a bounded domain whose boundary consists of a finite number of 
pairwise disjoint, real-analytic simple closed curves. Let 
$A(D)$ be the algebra of all continuous 
functions on $\overline D$ which are holomorphic on $D$. Denote by $H^1(D)$ the space 
of all holomorphic functions on $D$ such that $z\mapsto |h(z)|\ (z\in D)$ has a
harmonic
majorant [R]. Every
$h\in H^1(D)$ has nontangential boundary values $h^\ast $ almost everywhere on 
$bD$, \ $h^\ast\in L^1(bD)$ and
$$
h(z)={1\over{2\pi i}}\int _{bD}{{h^\ast (\z )d\z }\over{\z -z}}\ \ (z\in bD).
$$

We say that $f\in L^1(bD)$ extends holomorphically through $D$ if there
is $h\in H^1(D)$ such that 
$h^\ast =f $ almost everywhere on $bD$. 
 We say that $f\in L^1 (bD)$ extends meromorphically  through $D$ if there are a 
function $h\in H^1(D)$ and 
a nonzero polynomial $Q$ with all zeros contained in $D$ such that $f=h^\ast /Q$
almost 
everywhere on $bD$, 
or, equivalently, if $Qf$ extends holomorphically through $D$.

A function $f\in L^1 (bD)$ extends holomorphically through $D$ if and only if
$$
\int_{bD}f(\z )\omega (\z ) d\z = 0
$$
for each $\omega\in A(D)$ [R], which, since rational functions with 
poles outside $\overline D$ are dense in $A(D)$, is equivalent to 
$$
{1\over{2\pi i}}\int_{bD}{{f(\z )d\z }\over{\z -z}} \equiv 0\ \ (z\in
\C\setminus\overline D).
$$
There is no such simple test for meromorphic extendibility. If we happen to know the potential 
poles and their multiplicities, that is, if we know $Q$ then to check the meromorphic extendibility 
is easy as we simply 
check whether $Qf$ extends 
holomorphically through $D$. The problem becomes more difficult if we do not 
know in advance where the poles are. 
In this paper we show that if we know an upper bound for the number of poles then all possible 
candidates for $Q$ can be easily determined in advance so 
the question about meromorphic extendibility can be easily reduced to the 
question about holomorphic extendibility. 
\vskip 4mm
\bf 2.\ Poles of meromorphic extensions
\vskip 2mm \rm
Let $f\in L^1(bD)$. For large $z$ we have 
$$
{1\over{2\pi i}}\int_{bD}{{f(\z )d\z }\over{\z -z}}=
{{c_1(f)}\over{z}}+{{c_2(f)}\over{z^2}}+\cdots 
$$
where
$$
c_j(f)=-{1\over{2\pi i}}\int_{bD}\z^{j-1}f(\z )d\z \ \ (j\in \N).
\eqno (2.1)
$$
Given $N\in\N$ there is a nonzero polynomial $P(z)= D_0+D_1(z)+\cdots +D_Nz^N
$
such that for large $z$ we have
$$
P(z) {1\over{2\pi i}}\int_{bD}{{f(\z )d\z }\over{\z -z}}= 
\Biggl[ {{d_{N+1}}\over{z^{N+1}}}+ 
{{d_{N+2}}\over{z^{N+2}}} +\cdots\Biggr]  +R(z)
$$
where $R$ is a polynomial and $d_j,\ j\geq N+1$, are complex numbers. To get such
$P$ we have to solve the system
$$
\left.\eqalign{
&c_1(f)D_0+c_2(f)D_1+\cdots c_{N+1}(f)D_N = 0\cr
&\cdots \cr
&c_N(f)D_0+c_{N+1}(f)D_1+\cdots c_{2N}(f)D_N = 0 . \cr}\right\} 
\eqno (2.2)
$$
This is a homogeneous system of $N$ linear equations with $N+1$ unknowns 
which always has a nontrivial solution.
\vskip 2mm
\noindent\bf THEOREM 2.1\ \it Let $f\in L^1(bD)$ and let $c_j (f),\ j\in \N$, 
be as in \hbox{\rm (2.1)}. Let $N\in \N$ and let 
$P(z) = D_0+D_1 z +\cdots +D_N$ where $D_0, D_1,\cdots D_N$ is a nontrivial 
solution of the system \hbox{\rm (2.2)}.  
 The function $f$ extends meromorphically through $D$ with the extension 
 having at most $N$ poles in $D$, 
counting multiplicity, if and only if the function $z\mapsto P(z)f(z)$ 
extends holomorphically through $D$.
\vskip 2mm
\rm \noindent Thus, if we are asking whether $f$ has a meromorphic extension 
through $D$ with at most $N$ poles in $D$, 
counting multiplicity, then the only candidates for the poles are the zeros 
of $P$. For an analogous result for continuous 
functions on the unit circle see [G1]. 

Suppose that $z\mapsto P(z)f(z)\ (z\in bD)$ extends holomorphically through $D$. 
If $P$ has no zero on $bD$ then $f$ extends meromorphically through 
$D$ for in this case we can write $P=QS$ where 
the polynomial $Q$ has all its zeros on $D$ and the polynomial $S$ 
has all its zeros in $\C\setminus \overline D$. 
So if $Pf=h^\ast $ on $bD$ where $h\in H^1(D)$ then 
$$
f= {{h^\ast /S}\over Q} \hbox{\ \ almost everywhere on \ \ } bD
$$
where $h/S\in H^1(D)$. Even in the case when $P$ has zeros on $bD$ the function $f$
extends meromorphically through 
$D$. This follows from the following
\vskip 2mm
\noindent \bf LEMMA 2.2\ \it Let $f\in L^1(bD)$,\ let $a\in bD$ and assume
that $z\mapsto (z-a)f(z)\  (z\in bD)$ extends holomorphically through $D$. Then 
$f$ extends holomorphically 
through $D$.
\vskip 2mm
\noindent \rm In other words, if $g\in H^1(D)$ and if the function 
$z\mapsto g^\ast (z)/(z-a)\ (z\in bD)$ belongs to 
$L^1(bD)$ then there is a function $h\in H^1(bD)$ such that 
$h^\ast (z)= g^\ast (z)/(z-a)$ a.e.on $bD$.
\vskip 4mm
\bf 3.\ Some facts about the Cauchy integrals\rm

\vskip 2mm Our $bD$ consists of $m$ pairwise 
disjoint curves $\Gamma_1,\cdots ,\Gamma_m$ where $\Gamma _m$ is the
boundary of 
the unbounded component $D_m$ of $\C\setminus\overline D$. For each $j,\ 1\leq j\leq
m-1,$ 
let $D_j$ be the domain bounded by $\Gamma_j$. 

For each $j,\ 1\leq j\leq m-1$, let $G_j\in A(D_j)$, let $G\in A(D)$ and let $G_m$ 
be a continuous function on $\overline{D_m}$, holomorphic on $D_m$ and vanishing at
infinity. Assume that
$G$ and $G_j,\ 1\leq j\leq m$ all have smooth boundary values so that the function
$g$, defined on $bD$ by
$$
g(z)= G(z)-G_j(z)\ \ (1\leq j\leq m)
\eqno (3.1)
$$
is smooth. The function $g$ determines the functions $G$ and $ G_j,\ 1\leq j\leq m,$
uniquely. 
For, if $H_j\in A(D_j),\ 1\leq j\leq m-1$ and $H\in A(D)$ are functions with smooth
boundary 
values  and $H_m$ is smooth on $\overline{D_m}$, holomorphic on $D_m$ and vanishing
at infinity such that
$$
g(z)=H(z)-H_j(z)\ (1\leq j\leq m)
\eqno (3.2)
$$
then the function 
$$
\Phi (z) =\left\{ \eqalign{
& G(z)-H(z)\ \ (z\in\overline D)\cr
& G_j(z)-H_j(z)\ \ (z\in\overline {D_j},\ 1\leq j\leq m)\cr}\right.
$$
is, by (3.1) and (3.2), well defined, continuous on $\C$ and holomorphic 
on $\C\setminus bD$ and vanishing at infinity. So $\Phi $ is an entire function
vanishing at infinity 
hence $\Phi\equiv 0$ so $G\equiv H$ and $G_j\equiv H_j\ (1\leq j\leq m)$. In fact,
by the Plemelj jump formulas for the Cauchy integrals we have
$$
{1\over{2\pi i}}\int_{bD}{{g(\z )d\z }\over{\z -z}}=
\left\{ \eqalign{
& G(z)\ \ \ (z\in D)\cr
& G_j(z)\ \ (z\in D_j,\ 1\leq j\leq m).\cr}\right.
$$

Suppose now that $g(z) = R(z)/S(z)\ (z\in bD)$ 
where S is a polynomial with all zeros contained in $D$ and $R$ is a polynomial, $deg R < deg S$. Then the functions $$
G(z) = 0\ \ (z\in \overline D)
$$
$$
G_j(z) = -R(z)/S(z)\ \ (z\in \overline {D_j},\ 1\leq j\leq m)
$$ have the properties above so
$$
{1\over{2\pi i}}\int_{bD}{{g(\z )d\z }\over{\z -z}} = -R(z)/S(z)\ \ (z\in D_m).
$$

\bf 4.\ Proof of Theorem 2.1
\vskip 2mm
\rm We first prove Lemma 2.2.

The assumption implies that for each $z\in\C\setminus\overline D$ we have 
$$
0 = {1\over{2\pi i}}\int_{bD}{{(\z -a)f(\z )}\over{\z - z}} d\z = 
{1\over{2\pi i}}\int_{bD} f(\z )d\z + 
(z-a){1\over{2\pi i}}\int_{bD}{{f(\z )d\z }\over{\z - z}}
$$
hence
$$
{{1}\over{2\pi i}}\int_{bD}{{f(\z )d\z }\over{\z - z}} = - 
{1\over{z-a}}\Biggl[{1\over{2\pi i}}\int_{bD} f(\z )d\z 
\Biggr]\ \ (z\in \C\setminus  \overline D).
$$
We shall show that this implies that $\int _{bD}f(\z ) d\z = 0$ so
$$
{1\over{2\pi i}}\int_{bD}{{f(\z )d\z }\over{\z - z}} = 0\ \ 
(z\in \C\setminus\overline D)
\eqno (4.1)
$$
which we want to prove. Suppose, contrarily to what
we want to prove, that there is an $L\not= 0$ such that 
$$
{1\over{2\pi i}}\int_{bD}{{f(\z )d\z }\over{\z - z}} = 
{L\over{z-a}}\ \ (z\in \C\setminus \overline D).
$$
It follows that 
$$ \lim_{z\rightarrow a, z\in \C\setminus\overline 
D}(z-a) \int_{\Lambda}{{f(\z )d\z }\over{\z - z}} = L
\eqno (4.2)
$$
for each arc $\Lambda \subset bD$ containing $a$ in its interior. 
Denote by $\D $ the open unit disc. 
Since $bD$ is real analytic there are a map $\Phi $ mapping a disc 
$\Omega $ centered at the origin 
biholomorphically onto a neighbourhood of $a$, and $C,\ 0<C<\infty$, such that $\Phi (0)=a$, 
such that $\Phi$ maps $[-T,T]\subset \Omega$ 
onto an arc $\Lambda \subset bD$, the upper half of $T\D $ 
to $\C\setminus \overline D$, the lower half of $T\D$ to $D$, and 
such that
$$
{1\over C}|z-w|\leq |\Phi (z)-\Phi (w)|\leq C|z-w|\ \ (z,w\in\Omega ).
$$
Now, (4.2)implies that 
$$
\lim_{t\searrow 0}[\Phi (it)-\Phi (0)] \int_{-T}^T{{f(\Phi(x))\Phi^\prime 
(x)dx}\over{\Phi (x)-\Phi (it)}}\not= 0.
\eqno (4.3)
$$
Write $g(x) = f(\Phi (x))\Phi^\prime (x)$. Let $\varepsilon >0$. Since $f$ is integrable it follows that 
$g$ is
integrable so 
there is an $M<\infty $ such that 
$$
\int_{\{ x\colon |g(x)|\geq M\} } |g(x)| dx <\varepsilon .
$$
Let
$$ A_M = \{ x\in [-T,T]\colon \ |g(x)|<M\}\ \ \ \ \ B_M = \{ x\in [-T,T]\colon\ 
|g(x)|\geq M\} .
$$
Since $\Phi^\prime (0)\not= 0$,  (4.3) implies that
$$
\lim_{t\searrow 0} t\int_{-T}^T{{g(x)}\over{x-it}}\ {{x-it}\over{\Phi (x)-\Phi (it)}}
dx\not= 0
\eqno (4.4)
$$
We have
$$
\Biggl\vert t \int_{B_M}{{g(x)}\over{x-it}}\ {{x-it}\over{\Phi (x)-\Phi
(it)}}dx\Biggr\vert \leq |t|.{1\over{|t|}}.C.
\int_{B_M}|g(x)|dx\leq C\varepsilon .
$$
Further,
$$ 
\eqalign{
&\Biggl\vert t \int_{A_M}{{g(x)}\over{x-it}} . {{x-it}\over{\Phi (x)-\Phi (it)}} 
dx\Biggr\vert =
\Biggl\vert t \int_{A_M}g(x){{x+it}\over{x^2+t^2}}.{{x-it}\over{\Phi (x)-\Phi (it)}}
dx\Biggr\vert \cr
&\leq |t|.\int_{-T}^T M.C.{{|x|}\over{x^2+t^2}} dx +
|t|^2\int_{-T}^T{{MCdx}\over{x^2+t^2}} .\cr}
$$
It is easy to see that both terms the last expression tend to zero as $t\searrow 0$ so 
$$
\Biggr\vert \int_{-T}^T{{g(x)}\over{\Phi (x)-\Phi (it)}} dx\Biggr\vert \leq
2C\varepsilon
$$
provided that $t>0$ is small enough which contradicts (4.4) since $\varepsilon $ 
can be chosen arbitrarily small. This completes the proof of Lemma 2.2.

We now turn to the proof of the theorem. Let $P$ be as in Theorem 2.1 and 
assume that $Pf$ extends holomorphically through $D$. If $a\in bD$
is a zero of $P$ then by Lemma 2.2 the function $z\mapsto [P(z)/(z-a)]f(z)$ 
extends holomorphically through $D$.  
Thus, factoring out zeros of $P$ contained in $bD$ we conclude that there is 
a polynomial $Q$ having no zero on $bD$ and having at most $N$ zeros in $D$ 
such that 
$z\mapsto Q(z)f(z)$ extends holomorphically through $D$ so the function $f$ 
extends meromorphically through $D$ and the extension has 
at most $N$ poles in $D$, counting multiplicity.  This completes the 
proof of the if part. 

To prove the only if part of Theorem 2.1 observe first that if $h\in H^1(D)$ and
$a\in D$ then
$$
{{h(z)}\over{(z-a)^m}} = {{h(a)}\over{(z-a)^{m-1}}}+\cdots + 
{{h^{(m-1)}(a)}\over{(m-1)!(z-a)}} + w(z)
$$
where $w\in H^1(D)$. Assume now that $f$ extends meromorphically through $D$ 
with the extension having at most $N$ poles in $D$, counting multiplicity. 
This means that either $f$ extends holomorphically 
through $D$ - there is nothing to prove in this case - or there are 
$a_1,  \cdots a_J\in D$, and $k_1,\cdots k_J\in\N$, such that $k_1+\cdots k_J\leq N$ and that 
$$
f(z) ={{G(z)}\over{(z-a_1)^{k_1}\cdots (z-a_J)^{k_J}}}\ \ \ (z\in bD)
$$
where $G\in H^1(D)$.  Decomposing 
$$
1\over{(z-a_1)^{k_1}\cdots (z-a_J)^{k_J}}
$$
into partial fractions the 
preceding discussion implies that  
$$f(z)=H^\ast (z)+ T(z)\ \ (z\in bD)$$
where $H\in H^1 (D)$ and
$$
T(z) = {{R(z)}\over{S(z)}}\ \ (z\in bD)
$$
where $R, S$ are polynomials with no common zeros, $deg R< degS \leq N$ 
and all zeros of $S$ are contained in $D$.  Let $P$ be as
in Theorem 1.1. 
 We want to prove that $Pf = PH^\ast +PT $ extends holomorphically through $D$. 
Since $PH\in H^1 (D)$ it is enough to prove that $PT$ extends holomorphically through $D$.

From Section 3 we see that
$$
{1\over{2\pi i}}\int_{bD}{{f(\z )d\z }\over{\z - z}} =
{1\over{2\pi i}}\int_{bD}{{H^\ast (\z )d\z }\over{\z - z}} +
{1\over{2\pi i}}\int_{bD}{{T(\z )d\z }\over{\z - z}} = 0 - R(z)/S(z) \ (z\in D_m)
$$
which implies that for large $z$ the polynomial $P$ satisfies
$$
P(z)\biggl(-{{R(z)}\over{S(z)}}\biggr) = Q(z) + {{c_{N+1}}\over {z^{N+1}}}+{{c_{N+2}}\over {z^{N+2}}}+\cdots 
$$
where $Q$ is a polynomial. There are arbitrarily small $\alpha\in \C$ such that  
$$
P(z)\bigl({{R(z)}\over{S(z)}}\bigr) + Q(z)+\alpha \not= 0 \ \ \ (z\in bD).
\eqno (4.5)
$$ 
Assume for a moment that not all $c_j,\ j\geq N+1$, vanish. Then the left side of (4.1) equals
$$
\alpha - {{c_{N+1}}\over {z^{N+1}}}-\cdots 
$$
and by the argument principle the change of argument along $bD_m$ oriented as part of $bD$ 
is less that or equal to $-2\pi (N+1)$ 
provided that $\alpha$  is small enough. Thus, the change of argument of the left side of (4.5) along the part
of $bD$ that coincides with $bD_m$ is less that or equal to $-(N+1)2\pi$. Further, the left side of (4.5) is holomorphic 
on each $D_j,\ 1\leq j\leq m-1$, so by the argument principle its change of argument along the part of $bD$ that coincides 
with $bD_j$ is nonpositive, $1\leq j\leq m-1$. Thus, the change of argument of the left side of (4.5) along  $bD$ is less than or eaual to $- (N+1)2\pi$ 
which contradicts the fact that $S$ has at most $N$ zeros. This proves that $c_j=0\ (j\geq N+1)$ and so 
$$
P(z)T(z) = -Q(z)\ (z\in bD)
$$
so 
$PT$ extends holomorphically through $D$. The proof is complete. 
\vskip 4mm
\bf 5.\ Meromorphic extendibility of continuous functions\rm
\vskip 2mm

We show that for continuous functions on $bD$ the meromorphic 
extendibility can be expressed in terms of the argument principle. 

Given a continuous function  $\varphi\colon\ bD\rightarrow \C\setminus \{0\}$ we
denote by $\cW (\varphi )$ the winding number of 
$\varphi $ (around the origin). So $2\pi \cW (\varphi ) $ equals the
change of argument of $\varphi (z)$ as $z$ 
runs along $bD$ following the standard orientation. 

If a continuous function $f\colon\ bD\rightarrow \C\setminus \{ 0\} $ extends
meromorphically through $D$ 
then $\cW (f)\geq -N$ where $N$ is the number
of poles of the meromorphic extension $\tilde f$ (counted with multiplicity). 
Indeed, by the argument principle, 
$$
\cW (f) = \nu _0(\tilde f)-\nu _p(\tilde f) = \nu _0 (\tilde f) - N\geq -N
$$
where $\nu _0 (\tilde f)$ is the number of zeros of $\tilde f$ on $D $ 
and $\nu _p(\tilde f)$ is the number of poles of 
$\tilde f$ on $D$.

Let $f \colon \ bD\rightarrow \C$ be a continuous function which extends
meromorphically
through $D $ and whose meromorphic extension $\tilde f$ has 
$N$ poles on $D $. Then
$
\cW (Pf+Q) \geq -N 
$
for all functions $P, Q $ in $A(D)$  such that $Pf+Q\not = 0$ on $bD$. Indeed,
$P\tilde f + Q$, the 
meromorphic extension of $Pf+Q$, has no other poles than $\tilde f$ and therefore, 
by the argument principle, \ $\cW (Pf+Q)\geq -N$.  This property characterizes meromorphic extendibility: 
\vskip 2mm
\noindent \bf THEOREM 5.1\ \ \it Let $D$ be a bounded domain in $\C$ whose boundary consists of a finite number 
of pairwise disjoint simple closed curves. A continuous function $f\colon \ bD\rightarrow \C$ 
extends meromorphically through 
$D $ if and only if there is an $N\in\N\cup \{0 \} $ such that 
$$
\cW (Pf+Q)\geq -N
\eqno (5.1)
$$
for all  $P, Q \in A(D) $ such that $Pf+Q\not= 0 $ on $bD $. If this 
is the case then the meromorphic extension of $f$ has at most $N$ poles 
in $D$, counting multiplicity. \rm
\vskip 2mm
\noindent If $N=0$ there are no poles so  
we have holomorphic extendibility - when $D$ is the unit disc a better theorem (with $P\equiv 1$)  was proved by the author 
in [G4] and a simpler proof was obtained by D.\ Khavinson in [K], and if $D$ is a multiply connected domain, 
such a theorem 
was proved by the author in [G2]. 
If $D$ is the open unit disc the theorem was proved for general $N$ by the author in [G1] using 
 harmonic analysis. 
\vskip 2mm
\noindent\bf Proof.\ \rm The only if part follows from the argument 
principle as we have shown above. To prove the if part, assume that (5.1) holds for all 
$P, Q\in A(D)$ such that $Pf+Q\not=0$ on $bD$. Passing to a smaller $N$ if necessary we may assume that 
there are $P_0, Q_0\in A(D)$ such that $P_0 f+ Q_0\not=0$ on $bD$ and such that $\cW(P_0 f + Q)=-N$.
Write $P_0f+Q_0 = F$. Obviously $F\not=0$ on $bD$. Given $P, Q\in A(D)$ we have $PF+Q = P(P_0f+Q_0)+Q =
PP_0 f + (PQ_0+Q) $ where $PP_0\in A(D)$ and $PQ_0+Q\in A(D)$ so (5.1) implies 
that $\cW (PF+Q)\geq -N$ whenever $P, Q\in A(M)$ are such that $PF+Q\not=0$ on $bD$. It follows that 
$$
\cW \biggl(P+Q {1\over F}\biggr) =
 \cW \biggl({{PF+Q}\over F}\biggr) = \cW (PF+Q)-\cW (F) \geq -N - (-N) = 0
$$ 
whenever 
$P, Q\in A(D)$ are such that $P+Q{1\over F}\not= 0$ on $bM$. The main result of 
[G2] implies that $1/F$ extends holomorphically through $D$,
that is, $1/F = G$ where $G\in A(D)$. Obviously $G$ has no zero on $bD$ 
and since $\cW (F)=-N$ it follows that
$G$ has precisely $N$ zeros on $D$ counting multiplicity. Thus $1/G$ is a meromorphic extension of $F$ through $D$ 
which has no zero on $D$ and has exactly $N$ 
poles on $D$, counting multiplicity.  If $\alpha$ is suitably chosen small number then by the 
same process we see that $F_1 = 
(P_0 +\alpha)f+Q_0$ has also a meromorphic extension through $D$ with exactly $N$ poles in $D$. 
Thus, $f= (F_1-F)/\alpha $ extends 
meromorphically through $D$. We can choose $\alpha $ so that $P+\alpha\not=0$ at the poles of 
the meromorphic extension of $f$ so this meromorphic extension has the same number of poles as the meromorphic 
extension of $F_1$ has the 
same number of poles as $F_1$ which is $N$.
\vskip 2mm
For holomorphic extendibility (N=0) this theorem holds also for open Riemann surfaces [G3] 
and therefore, repeating the proof, we see
that Theorem 5.1 holds for open Riemann surfaces as well.

 We have already mentioned that there is a better version of Theorem 5.1 in the case when $N=0$:\it \ A continuous function $f$ on $bD$ 
extends holomorphically through $D$ if and only if \ $\cW (f+Q)\geq 0$ for all $Q\in A(D)$ such that 
$f+Q\not= 0$ on $bD$ \rm [G2]. That is, in the case when $N=0$ one can take $P\equiv 1$. 
Whether this holds for $N\geq 1$ remains to be seen:
\vskip 2mm
\noindent\bf Question \rm Let $N\in\N$ and suppose that $f$ is a 
continuous function on $bD$ such that 
$\cW (f+Q)\geq -N$ whenever $Q\in A(D)$ is such that $f+Q\not= 0$ on $bD$. 
Does it follow that $f$ extends meromorphically through $D$?
\vskip 2mm
\noindent The answer is not known even in the case when $D$ is a disc.

\vskip 5mm
This work was supported 
in part by the Ministry of Higher Education, Science and Technology of Slovenia 
through the research program Analysis and Geometry, Contract No.\ P1-0291
\vfill
\eject
\centerline{\bf REFERENCES}
\vskip 3mm

\noindent [F]\ S.\ D.\ Fisher:\ \it Function Theory on Planar Domains.

\noindent \rm  John Wiley \& Sons, New York 1983
\vskip 2mm
\noindent [G1]\ J.\ Globevnik:\ Meromorphic extendibility and the argument principle.

\noindent Publ.\ Mat.\ 52 (2008) 171-188
\vskip 2mm
\noindent [G2]\  J.\ Globevnik:\ The argument principle and holomorphic extendibility.

\noindent Journ.\ d'Analyse.\ Math.\ 94 (2004) 385-395
\vskip 2mm
\noindent [G3]\ J.\ Globevnik:\ The argument principle and holomorphic extendibility 
to finite Riemann surfaces.

\noindent  Math.\ Z.\ 253 (2006) 219-225
\vskip 2mm
\noindent [G4]\ J.\ Globevnik:\ Holomorphic extendibility and the argument principle.

\noindent Contemp.\ Math.\  Vol.\ 382 (2005) 171-175
\vskip 2mm
\noindent [K]\ D.\ Khavinson:\ A note on a theorem of J. Globevnik

\noindent Contemp.\ Math.\  Vol.\ 382 (2005) 227-228
\vskip 2mm

\noindent [R]\ W.\ Rudin:\ Analytic functions of class $H_p$.

\noindent Trans.\ Amer.\ Math.\ Soc.\ 78 (1955) 46-66
\vskip 2mm
\noindent [S]\ E.\ L.\ Stout: \it The theory of uniform algebras.\rm 

\noindent Bogden and Quigley, Tarrytown-on-Hudson, NY, 1971
\vskip 2mm

\noindent [Z] A.\ Zygmund: \it Trigonometric series. \rm

\noindent Cambridge University Press, Cambridge, New York, 1959

\vskip 8mm
\noindent Institute of Mathematics, Physics and Mechanics

\noindent University of Ljubljana, Ljubljana, Slovenia

\noindent josip.globevnik@fmf.uni-lj.si

\bye